\documentclass[multi]{cambridge7A}
\usepackage[UKenglish]{babel}
\usepackage[longnamesfirst,sectionbib]{natbib}
\usepackage{chapterbib}
\usepackage{amsmath,amssymb,amsthm,amsfonts}
\usepackage{enumerate,url}

\setcitestyle{author,year,round,semicolon}

\theoremstyle{plain}
\newtheorem {theorem}{Theorem}[section]
\newtheorem {corollary}[theorem]{Corollary}
\newtheorem {proposition}[theorem]{Proposition}

\allowdisplaybreaks[1]

\newcommand{{\Polya}}{P\'olya{} }

\def \a   {\alpha}

\def \M{{\cal M}}

\def \eps {\varepsilon}

\def\P{{\mathbb{P}}}
\def \R {{\mathbb{R}}}
\newcommand{\Z}{\mathbb{Z}}
\def\F{{\cal{F}}}

\def\limt{\lim_{t\to\infty}}
\def\limt0{\lim_{t\to 0}}

\def\|{\,|\,}
\def \eps {\varepsilon}

\def\bn{\begin{eqnarray*}}
\def\en{\end{eqnarray*}}
\def\bnn{\begin{eqnarray}}
\def\enn{\end{eqnarray}}

\newcommand{\ri}{\rightarrow}

\providecommand*\Index[1]{#1\index{#1}}
\providecommand*\undex[1]{} 

\begin {document}
\alphafootnotes
\author[V. Shcherbakov and S. Vol\-kov]{Vadim Shcherbakov\footnotemark\ and
  Stanislav Volkov\footnotemark }
\chapter[Queueing with neighbours]{Queueing with neighbours}
\footnotetext[1]{Laboratory of Large Random Systems, Faculty of Mechanics
  and Mathematics, Moscow State University, 119991 Moscow, Russia;
  v.shcherbakov@mech.math.msu.su}
\footnotetext[2]{Department of Mathematics, University of Bristol, University
  Walk, Bristol BS8~1TW; S.Volkov@bristol.ac.uk}
\arabicfootnotes
\contributor{Vadim Shcherbakov \affiliation{Moscow State University}}
\contributor{Stanislav Volkov \affiliation{University of Bristol}}
\renewcommand\thesection{\arabic{section}}
\numberwithin{equation}{section}
\renewcommand\theequation{\thesection.\arabic{equation}}
\numberwithin{table}{section}
\renewcommand\thetable{\thesection.\arabic{table}}

\begin {abstract}
In this paper we study asymptotic behaviour of a growth process
generated by a semi-deterministic variant of the cooperative
sequential adsorption  model (CSA). This model can also be viewed
as a particular example from queueing theory, to which John
Kingman has contributed so much. We show that the quite limited
randomness of our model still generates a rich collection of
possible limiting behaviours.
\end {abstract}

\subparagraph{Keywords}cooperative sequential adsorption,
interacting particle systems, max-plus algebra, queueing, Tetris

\subparagraph{AMS subject classification (MSC2010)}Primary 60G17, 62M30;
Secondary\break 60J20

\section{Introduction}

Let $\M=\{1,2,\ldots,M\}$ be a lattice segment  with periodic
boundary conditions (that is, $M+1$ will be understood as $1$ and
$1-1$ will be understood as $M$), where $M\geq 1$. The growth
process studied in this paper is defined as  a discrete-time
\index{Markov, A. A.!Markov chain|(}Markov chain $(\xi_i(t),\, i\in\M,\, t\in \Z_{+}),$ with values in
$\Z^{M}_{+}$ and specified by the following transition
probabilities:
 \bnn\label{tran}
 &\P\left(\xi_{i}(t+1)=\xi_{i}(t)+1,\ \xi_j(t+1) = \xi_j(t)\ \forall j\ne i\|\xi(t)\right) & \nonumber\\
 &=\left\{\begin{array}{rcl}
  0, &\mbox{if} & u_i(t)>m(t),\\
  1/N_{\min}(t), &\mbox{if} & u_i(t)=m(t),
 \end{array}\right. &
 \enn
for $i\in\M$, where
$$
u_i(t)=\sum\limits_{j\in U_i} \xi_j(t),\,\, i\in\M,
$$
 $U_i$ is
a certain neighbourhood of site $i$,
\begin{equation}
\label{mt}
m(t)=\min_{k\in\M}u_k(t)
\end{equation}
and  $N_{\min}(t)\in\{1,2,\dots,M\}$ is the number of $u_i(t)$
equal to $m(t)$. The quantity $u_i(t)$ is called  the
\index{potential|(}\emph{potential} of site $i$ at time $t$.

The growth process describes the following random sequential
allocation procedure. Arriving particles are sequentially
allocated at sites of $\M$ such that a particle is allocated
uniformly over sites with minimal potential. Then the process
component $\xi_k(t)$ is the number of particles at site $k$ at
time $t$. The growth process can be viewed as a certain limit case
of a growth process studied  in~\cite{ShchVolk}. The growth
process in~\cite{ShchVolk} is defined as a discrete-time Markov
chain $(\xi_i(t),\, i\in\M,\, t\in \Z_{+}),$ with values in
$\Z^{M}_{+}$ and specified by the following transition
probabilities
\begin{equation}
\label{growth}
\P\{\xi_i(t+1)=\xi_{i}(t)+1,\ \xi_j(t+1)=\xi_j(t)\ \forall j\ne
i\|\xi(t)\}= \frac{\beta^{u_i(t)}}{Z(\xi(t))}
\end{equation}
where
$$Z(\xi(t))=\sum_{j=1}^{M}\beta^{u_j(t)},$$
$\beta$ is a positive  number  and the other notations are the same as
before. It is easy to see that the process defined by transition
probabilities (\ref{tran}) is  the corresponding limit process  as
$\beta\ri 0$. In turn, the growth process specified by the
transition probabilities (\ref{growth}) is a particular version of
the
\Index{cooperative sequential adsorption model (CSA)}. CSA is a
probabilistic model which is widely used in physics for modelling
various
\index{adsorption|(}adsorption processes (see
\index{Evans, J. W.}\cite{Evans},
\index{Privman, V.}\cite{Privman} and
references therein).  Some asymptotic and statistical studies  of
similar CSA in a continuous setting were undertaken in \cite{Sch},
\index{Penrose, M. D.}\cite{PenSch} and \cite{PenSch1}.

In~\cite{ShchVolk} we consider the following variants of
neighbourhood
 \bn
&{\bf (A1)}\quad {}& U_i=\{i\}          \mbox{ (empty),}\\
&{\bf (A2)}\quad {}& U_i=\{i, i+1\}     \mbox{ (asymmetric),}\\
&{\bf (A3)}\quad {}& U_i=\{i-1, i, i+1\} \mbox{ (symmetric),}
 \en
where, due to the periodic boundary conditions, $U_{M}=\{M,1\}$ in
case  {\bf (A2)} and $U_{1}=\{M, 1, 2\},\,U_{M}=\{M-1, M, 1\}$ in
case {\bf (A3)} respectively. It is easy to see that for the
growth  process  studied in this paper the case {\bf (A1)} is
trivial. Therefore,   we will consider cases {\bf (A2)} and {\bf
(A3)} only.

A stochastic process $u(t)=(u_1(t),\dots,u_{M}(t))$ formed by the
sites' potentials plays an important role in our asymptotic study
of the growth process. It is easy to see that $u(t)$ is also a
Markov chain, with transition probabilities given by
\begin{equation}\label{eq_tran}
\P\left(u_{i}(t+1)=u_{i}(t)+1_{\{i\in U_k\}}\right)=
  \begin{cases}0,&\text{if $u_k(t)>m(t)$},\\
     1/N_{\min}(t),&\text{if $u_k(t)=m(t)$},
  \end{cases}
\end{equation}
for $k\in\M$. This process has the following rather obvious
\index{queue}queueing interpretation
\index{Foss, S. G.}(S.~Foss, personal communications)
explaining  the title of the paper (originally titled `Random
sequential adsorption at extremes'). Namely, consider  a system
with $M$ servers, with the clients arriving in bunches of
 $2$  in case {\bf (A2)} and of  $3$ in case {\bf (A3)}. The quantity
$u_i(t)$ is interpreted as the number of clients at server $i$ at
time $t$.
In case {\bf (A2)}, of the two clients in the arriving pair, one joins the
shortest queue, the other its left neighbouring queue, the two choices being
equally likely. In case {\bf (A3)}, of the three clients in an arriving
triple, one joins the shortest queue, the others its left and right
neighbouring queues, the choices being equally likely.

Our goal is to describe the long time behaviour of the growth
process, or, equivalently, to describe the limiting queueing
profile of the network. It should be noted that the method  of
proof in this paper is purely combinatorial. This is in contrast
with ~\cite{ShchVolk}, where the results  are proved by combining
the
\Index{martingale} techniques from
\index{Fayolle, G.}\index{Malyshev, V. A.}\index{Menshiikov, M. V.}\cite{fmm}
with some probabilistic techniques used in the theory of
\index{reinforced process}reinforced processes from
\index{Tarres, P.@Tarr\`es, P.}\cite{Tarres} and~\cite{V2001}.

Observe that the model considered here can be viewed as a
randomized
\index{Tetris@\emph{Tetris}}\emph{Tetris} game, and hence it can possibly be
analyzed using the techniques of
\index{maxplus algebra@max-plus algebra}max-plus algebra as well; see
\index{Bousch, T.}\index{Mairesse, J.}\cite{BM} and Section 1.3 of
\index{Heidergott, B.}\index{Oldser, G. J.}\index{Woude, J. van der}\cite{HOW}
for details.

For the sake of completeness, let us mention  another limit case
of the growth process  specified by transition probabilities
(\ref{growth}): namely,  the limit process arising as $\beta\ri
\infty$. It is easy to see that the limit process in this case
describes the
\Index{allocation process} in which a particle is  allocated
with equal probabilities to one of the sites with maximal potential. The
asymptotic behaviour (as $t\ri \infty$) of this limit process  is
briefly discussed in Section \ref{binf}.

\section{Results}\label{res}

\begin{theorem}\label{thm:A2_b0}
Suppose $U_i=\{i,i+1\}$, $i\in\M$. Then, with
probability $1$, there is a $t_0=t_0(\omega)$ (depending also on
the initial configuration) such that for all $t\ge t_0$
 \bnn\label{eq_at_most2}
 |\xi_i(t)-\xi_{i+2}(t)|\le 2
 \enn
for $i\in\M$. Moreover,
 \bnn\label{eq_xieo}
 \xi_i(t)=\frac t{M} + \eta_i(t)+\left\{\begin{array}{rc}
  0, &\mbox{when $M$ is odd},\\ \\
   (-1)^i \, Z(t), &\mbox{when $M$ is even},
 \end{array}\right.
 \enn
where $|\eta_i(t)|\le 2M$ and for some $\sigma>0$
$$
\lim_{n\to\infty} \frac{Z(\lfloor sn\rfloor)}{\sigma
\sqrt{n}}=B(s),
$$
where $\lfloor x\rfloor$ denotes the integer part of $x$ and
$B(s)$ is a standard
\index{Brown, R.!Brownian motion}Brownian motion.
\end{theorem}

\begin{theorem}\label{thm:A3_b0}
Suppose $U_i=\{i-1,i,i+1\}$, $i\in\M$. Then with
\undex{almost sure convergence@almost-sure convergence}probability $1$
there exists the limit ${\bf x}=\lim_{t\ri \infty}\xi(t)/t$, which
takes a finite number of possible values with positive
probabilities. The set of limiting configurations consists of
those ${\bf x}=(x_1,\dots,x_{M})\in\R^M$ which simultaneously
satisfy the following properties:
\begin{itemize}
\item there exists an $\a>0$ such that $x_i\in\left\{0,\a/2,\a\right\}$ for all
$i\in\M$; also $\sum_{i=1}^{M} x_i=1$;
\item if $x_i=0$, then $x_{i-1}>0$ or $x_{i+1}>0$, or both;
\item if $x_i=\a/2$, then
$$(x_{j-3},x_{j-2},x_{j-1},x_{j},x_{j+1},x_{j+2})=(\a,0,\a/2,\a/2,0,
\a),$$
where $j\in\{i,i+1\}$;
\item if $x_i=\a$, then $x_{i-1}=x_{i+1}=0$;
\item if $M=3K$ is divisible by $3$, then
$$\min\{ x_{j},x_{j+3},x_{j+6},\dots,x_{j+3(K-1)}\}=0,$$
 for $j=1$, $2$, $3$.
\end{itemize}
Moreover, the
\index{adsorption|)}adsorption eventually stops at all $i\in\M$ where
$x_i=0$, that is
$$
 \sup_{t\ge 0}\xi_i(t)=\infty \mbox{ if and only if } x_i>0.
$$
Additionally, if  the initial configuration is empty, then for
each $x_i=0$ we must have that {\bf both} $x_{i-1}>0$ {\bf and}
$x_{i+1}>0$.
\end{theorem}

\begin{table}[htb]
\caption{Limiting configurations for symmetric interaction}
\begin{center}
\begin{tabular}{|c|c|c|}
\hline
  M  & Limiting configurations (up to rotation)& No. of limits \\
 \hline
 $4$   & {\small$\left(\frac 12,0,\frac 12,0\right)$} & 2  \\
 \hline
 $5$   &
 {\small
 $\left(\frac 14,\frac 14,0,\frac 12,0\right)$, $\left(\frac 12,0,0,\frac 12,0\right)^{*}$ } & 5  (10${}^{*}$) \\
 \hline
 $6$   & {\small$\left(\frac 13,0,\frac 13,0,\frac 13,0\right)$ } & 2  \\
 \hline
 $7$   &
 {\small
    $\left(\frac 16,\frac16,0,\frac 13,0,\frac 13,0\right)$,
         $\left(\frac 13,0,\frac13,0,\frac 13,0,0\right)^{*}$}
 & 7  (14${}^{*}$)\\
 \hline
 $8$   &
 {\small
  $\left(\frac 14,0,\frac14,0,\frac 14,0,\frac  14,0\right)$,} & 2  (18${}^{*}$) \\
  &
  {\small
   $\left(\frac 16,\frac 16,0,\frac13,0,0,\frac 13,0\right)^{*}$,
    $\left(0,0,\frac 13,0,0,\frac 13,0,\frac13\right)^{*}$}
  & \\
 \hline
 $9$   &
 {\small
 $\left(\frac 18,\frac18,0,\frac 14,0,\frac 14,0,\frac 14,0\right)$,
       $\left(0,0,\frac 14,0,\frac 14,0,\frac 14,0,\frac14\right)^{*}$
       }
       & 9 (18${}^{*}$) \\
 \hline
 $10$   &
 {\small  $\left(\frac 18,\frac18,0,\frac 14,0,\frac 18,\frac 18,0,\frac 14,0\right)$,
   $\left(\frac 15,0,\frac 15,0,\frac 15,0,\frac 15,0,\frac 15,0\right)$,}  &
   7 (42${}^{*}$)
   \\
   &
{\small
 $\left(0,\frac 14,0,0,\frac 14,0,\frac 14,0,\frac 18,\frac 18\right)^{*}$,
   $\left(0,\frac 14,0,\frac 14,0,0,\frac 14,0,\frac 18,\frac   14\right)^{*}$,}
   &
   \\
   &
{\small
    $\left(0,0,\frac 14,0,0,\frac 14,0,\frac 14,0,\frac    14\right)^{*}$,
    $\left(0,0,\frac 14,0,\frac 14,0,0,\frac 14,0,\frac 14\right)^{*}$
    }
  &    \\
 \hline
\end{tabular}
\end{center}
\label{Tabl2}
\end{table}
\begin{table}[htb]
\begin{minipage}{106mm}
\caption{Numbers of limiting configurations for symmetric
interaction for larger $M$}
\end{minipage}
\footnotesize
\begin{tabular}{@{}ccccccc@{}}
\hline
  M&11&12&13&14&15&16\\\hline
 Distinct conf.&1(4${}^{*}$)&2(7${}^{*}$)&1(8${}^{*}$)&3($12{}^{*}$)&2(16${}^{*}$)&3(20${}^{*}$)\\\hline
All conf.&11(44${}^{*}$)&14(74${}^{*}$)&13(104${}^{*}$)&23(142${}^{*}$)&20(220${}^{*}$)&34(290${}^{*}$)\\\hline
\end{tabular}
\label{Tabl3}
\end{table}

We will derive the asymptotic behaviour of the process $\xi(t)$ from
the asymptotic behaviour of the process of potentials.  In turn
the study of the process of potentials is greatly facilitated  by
analysis of the following auxiliary process
 \bnn\label{eqv(t)}
 v_k(t)=u_k(t)-m(t),\ k=1,\ldots,M.
 \enn
Observe that $v(t)$ also forms a Markov chain on $\{0,1,2,\dots\}$
and for each $t$ there is a $k$ such that $v_k(t)=0$. Loosely
speaking, the quantities $v_k(t)$, $k=1$, \ldots, $M$, represent what
happens `on the top of the growth  profile'.

It turns out that in the asymmetric case there is a single  class
of recurrent states to which the  chain eventually falls and then
stays in forever.   As we show later, this class is a certain
subset of the set $\{0,1,2\}^M$  containing the origin $(0,\ldots,0)$.
Thus a certain `stability' of  the process of potentials is
observed as time goes to infinity.

In particular, it yields, as we show,  that there will not be long
queues in the system if $M$ is odd; however, this does not prevent
occurrence of relatively long queues if $M$ is even. For instance,
if $M$ is even, then one can easily imagine the profile with peaks
at even sites, and zeros at odd sites. Besides, observe that  the
sum of the potentials of the  even sites equals the sum  of
the potentials of the odd sites (see Proposition \ref{prop:even});
therefore the difference between the total queue to the even
sites, and the total
\Index{queue} to the odd ones, behaves similarly to
the zero-mean
\index{random walk (RW)}random walk. It means that there are rather long
periods of time during which much longer queues are observed at
the even sites in comparison with the queues at the odd sites, and
vice versa. Thus, in the case of the asymmetric interaction we
observe in the limit $t\ri \infty$ a  `comb pattern' when $M$ is
even, and a `flat pattern' when $M$ is odd.

The picture is completely different in the symmetric case. The Markov
chain $v(t)$ is transient for any $M$; moreover,  there can be
only finitely many paths along which the chain escapes to
infinity. By this we mean that if the chain reaches a state
belonging to a particular \emph{escape path}, then it will never
leave it and will go  to infinity along this path, and we will
show that there can be only a finite number of limit patterns. An
escape path can be viewed as   \emph{an attractor}, since similar
effects are observed in
\index{network!neural network}neural network models studied in
\index{Karpelevich, F. I.}\index{Malyshev, V. A.}\index{Rybko, A. N.}\cite{KMR}
and
\index{Turova, T. S.}\cite{MalTur}. In fact, the Markov chain $v(t)$
describes the same dynamics as the neural network models in \cite{KMR} and
\cite{MalTur} though in a slightly different set-up. The
difference seems to be technical but it results in  quite
different model behaviour. We do not investigate this issue  in
depth here.

Table~\ref{Tabl2} contains the list of all possible limiting
configurations (for proportions of customers/particles)  for small
$M$, while in Table~\ref{Tabl3} we provide only the numbers of
possible limiting configurations for some larger $M$. Note that in
the tables the symbol ${}^*$ stands for the configurations which
cannot be achieved from the empty initial configuration.
Unfortunately, we cannot compute exact numbers of possible
limiting configurations for any $M$; nor can we predict which of
them will be more likely (though it is obvious that if we start
with the empty initial configuration, all possible limits which
can be obtained by a rotation of the same ${\bf x}$ will have the
same probability.)

\section{Asymmetric interaction}
In the asymmetric case the
\index{potential|)}potential of site $k$ at time $t$ is
$$
u_k(t)=\xi_k(t)+\xi_{k+1}(t),\, k\in\M.
$$
The transition probabilities of the Markov chain
$u(t)=(u_1(t),\dots,\allowbreak u_{M}(t))$ are given by
\begin{multline*}
 \P\left(u_{i}(t+1)=u_{i}(t)+1_{i\in\{k-1,k\}}, i=1,\ldots,M| u(t)\right)\\
=\left\{\begin{array}{rcl}
  0, &\mbox{if} & u_k(t)>m(t),\\
  N^{-1}_{\min}(t), &\mbox{if} & u_k(t)=m(t),
 \end{array}\right.
\end{multline*}
for $k\in\M$, where $N_{\min}(t)\in\{1,2,\dots,M\}$ is the number
of $u_i(t)$ equal to $m(t)$.
\begin{proposition}\label{prop:even}
If $M$ is odd, then for any $u=(u_1,u_2,\dots,u_M)$ the system
 \bnn\label{eq:uxiasym}
  u_1&=&\xi_1+\xi_{2}\nonumber\\
  u_2&=&\xi_2+\xi_{3}\nonumber\\
  &\vdots&\\
  u_M&=&\xi_M+\xi_{1}\nonumber
 \enn
has a unique solution. On the other hand, if $M$ is even,
system~(\ref{eq:uxiasym}) has a solution if and only if
 \bnn\label{eq:evenodd}
 u_1+u_3+\dots+u_{M-1}= u_2+u_4+\dots+u_{M}.
 \enn
\end{proposition}

\begin{proof}For a fixed $\xi_1$ we can express the remaining
$\xi_k$ as
 \bnn\label{eq:soluxi}
  \xi_k&=&u_{k-1}-u_{k-2}+u_{k-3}-\dots+(-1)^{k-1} \xi_1,
 \enn
 for any $k=2$, \dots, $M$.
Now, when $M$ is odd, there is a unique choice of
 \bn
  \xi_1&=&\frac   12\left(u_M-u_{M-1}+u_{M-2}+\dots-u_2+u_1\right).
 \en
When $M$ is even, by summing separately odd and even lines of
(\ref{eq:uxiasym}) we obtain condition~(\ref{eq:evenodd}). Then
it turns out that we can set $\xi_1$ to be any real number, with
$\xi_k$, $k\ge 2$, given by (\ref{eq:soluxi}).
\end{proof}

The
following statement immediately follows from
Proposition~\ref{prop:even}.
\begin{corollary}\label{cor_even}
If $M$ is even, then
 \begin{align}\label{condi}
v_1(t)+v_3(t)+\dots+v_{M-1}(t)=v_2(t)+v_4(t)+\dots+v_{M}(t) \ \ \
\forall t.
 \end{align}
\end{corollary}

In the following two Propositions we will show that when either
$M$ is odd \emph{or} $M$ is even  and condition (\ref{condi})
holds, the state ${\bf 0}=(0,\dots,0)$ is recurrent for the Markov
chain $v(t)$. First, define the following stopping times
 \bnn\label{eqt_j}
 t_0&=&0, \nonumber\\
 t_j&=&\min\{t>t_{j-1}:\ m(t_j)>m(t_{j-1})\}, \ j=1,2,\dots.
 \enn
Let
$$
S(j)=\sum_{k=1}^{M} v_k^*(t_j)
$$
where
 \bn
 v_k^*(t_j)=\left\{\begin{array}{rl}
  v_k(t_j), &\mbox{if }  v_k(t_j)\ge 2,\\
  0, &\mbox{otherwise,}
 \end{array}\right.
 \en
where the stopping times are defined by (\ref{eqt_j}).

\begin{proposition}\label{prop_a_men}
$S(j+1)\le S(j)$. Moreover, there is an integer $K=K(M)$ and an
$\eps>0$ such that
$$
\P(S(j+K)-S(j)\le -1\| v(t_j))\ge \eps
$$
on the event $S(j)>0$.
\end{proposition}

\begin{proof}For simplicity, let us write $v_k$ for $v_k(t_j)$.
Take some non-zero element $a\ge 1$ in the sequence of $v_k$ at
time $t_j$. Whenever it is followed by a consecutive chunk of
$0$s, namely
$$
\dots  a\ \underbrace{0\ 0\ \dots \ 0}_m \dots
$$
at time $t_{j+1}$ this becomes either
$$
 \dots  a\ \underbrace{z_1\ z_2\ \dots \ z_m}_m\dots  $$
 or
$$  \dots  a-1\ \underbrace{z_1\ z_2\ \dots \ z_m}_m\dots, $$
where $z_i\in\{0,1\},$
and the latter occurs if the second $0$ from the left is chosen
before the first one. On the other hand, if $a$ is succeeded by a
non-zero element, say `$\dots\  a\ b\ \dots$' at time $t_{j+1}$
this becomes either `$\dots\ a-1 \ b \ \dots $' or `$\dots\ a-1
\ b-1 \ \dots $'. In all cases, this leads to $S(j+1)\le S(j)$.

Secondly, from the previous arguments we see that if there is at
least one $a\ge 2$ in the sequence of $(v_1 \dots v_{M})$ followed
by a non-zero element, then this element becomes $a-1$ at
$t_{j+1}$ and hence $S(j+1)\le S(j)-1$.

Now let us investigate what happens if the opposite occurs. Then
each element $a\ge 2$ in $(v_1 \dots v_{M})$ is followed by a
sequence of $0$s and $1s$ such that we observe either
$$ \dots  a\ 0\ b\dots $$
or
$$ \dots  a\ 0\ 1\ \underbrace{z_3\ z_4\ \dots \ z_m}_{m-2}\ b\dots, $$
where $m\ge 2$, $b\ge 2$ and $z_i\in\{0,1\}$.
Because of Corollary~\ref{cor_even}, we cannot have an alternating
sequence of $0$s and non-zero elements; therefore, we must be able
to find somewhere in the sequence of $v$s a chunk which looks
either like
\begin{equation}
 \dots  a\ \underbrace{\underbrace{0\ 1}\ \underbrace{0 \ 1}\ \dots\  \underbrace{0 \ 1}}_{2l\mbox{ elements}}\ 0\ 1\ c\dots
\mbox{ where } c\ge 1 \tag{A}
\end{equation}
or
\begin{equation}
 \dots  a\ \underbrace{\underbrace{0\ 1}\ \underbrace{0 \ 1}\ \dots\  \underbrace{0 \ 1}}_{2l\mbox{ elements}}\ \underbrace{0\ \ 0\ \ ?}_{\ i\ \ i+1\ i+2}\dots \tag{B}
\end{equation}
where $l\ge 0$. Note that a configuration of type (A) at time
$t_{j+1}$ with probability $1$ becomes a configuration of type
(B). At the same time, in configuration (B), with probability of
at least $\frac 13$ the $0$ located at position $i+1$ is chosen
before either the $0$ at position $i$ or (possibly) the $0$ at
position $i+2$. On this event, the configuration in (B) at time
$t_{j+1}$ becomes
\begin{equation}
 \dots  a\ \underbrace{\underbrace{0\ 1}\ \underbrace{0 \ 1}\ \dots\  \underbrace{0 \ 1}}_{2l-2\mbox{ elements}}\ 0 \ 0\ 0 \ ? \dots  \tag{$\mathrm{B}'$}
\end{equation}
By iterating this argument until $l=0$, we conclude that
eventually there will be  a chunk `$\dots a\ 0\ 0\dots$' on some
step $t_{j'}$ which in turn at time $t_{j'+1}$ will become `$
\dots a-1\ 0\ ?\dots $'with probability at least $\frac 13$,
resulting in $S_{j'+1}\le S_{j'}-1$. This yields the statement of
Proposition~\ref{prop_a_men} with $K=M$ and $\eps=3^{-M}$.
\end{proof}

\begin{proposition}\label{prop_a_conv}
With probability $1$, there is a $j_0=j_0(\omega)$ such that
$$
S(j)=0 \mbox{ for all } j\ge j_0.
$$
Additionally, the state ${\bf 0}=(0,0,\dots,0)$ is recurrent for the
Markov chain $v(t)$.
\end{proposition}

\begin{proof}The first statement trivially follows from
Proposition~\ref{prop_a_men}. Next observe that at times $t_j\ge
t_{j_0}$ the sequence $(v_1(t_j),\dots,v_{M}(t_j))$ consists only
of $0$s and $1$s locally looking either like
\begin{equation}
 \dots  1\ \underbrace{\underbrace{0\ 0}\ \underbrace{0 \ 0}\ \dots\  \underbrace{0 \ 0}}_{2l\mbox{ elements}}\ 1 \dots  \tag{C}
\end{equation}
or
\begin{equation}
 \dots  1\ \underbrace{\underbrace{0\ 0}\ \underbrace{0 \ 0}\ \dots\  \underbrace{0 \ 0}}_{2l\mbox{ elements}}\ 0\ 1 \dots  \tag{D}
\end{equation}
With positive probability  even-located $0$s are picked before
odd-located $0$s, hence at time $t_{j+1}$ configuration $(C)$
becomes
\begin{equation}
 \dots  0\ \underbrace{\underbrace{0\ 0}\ \underbrace{0 \ 0}\ \dots\  \underbrace{0 \ 0}}_{2l\mbox{ elements}}\ ? \dots  \tag{$\mathrm{C}'$}
\end{equation}
while configuration $(\mathrm{D})$ becomes
\begin{equation}
 \dots  0\ \underbrace{\underbrace{0\ 0}\ \underbrace{0 \ 0}\ \dots\  \underbrace{0 \ 0}}_{2l-2\mbox{ elements}}\ 0 \ 1\ 0\ ? \dots  \tag{$\mathrm{D}'$}
\end{equation}
In both cases $(\mathrm{C})\to(\mathrm{C}')$ and $(\mathrm{D})\to(\mathrm{D}')$ the number of $1$s
among the $v_i$ does not increase, and in the first case it goes
down by $1$. However, it is easy to see that whether $M$ is odd or
even (in the latter case due to Corollary~\ref{cor_even}) there
will be at least one chunk of type $(\mathrm{C})$, and hence with
positive probability $v(t)$ reaches state ${\bf 0}$ in at most
$M^2$ steps (since $t_{j+1}-t_j\le M$\,). The observation that
after $t_{j_0}$ the Markov chain $v(t)$ lives on a finite state
space $\{0,1,2\}^{M}$ finishes the proof.
\end{proof}

\begin{proof}[Proof of Theorem \ref{thm:A2_b0}]The first part easily
follows from Proposition~\ref{prop_a_conv} and the definition of
\index{potential|(}potentials $v$. Indeed, for $j\ge j_0$ and all $i$ we have
$v_i(t_j)\in\{0,1\}$, while for $t\in (t_j,t_{j+1})$ we have
$v_i(t)\in\{0,1,2\}$. On the other hand, omitting $(t)$, we can
write $v_{i+1}-v_{i}=u_{i+1}-u_i=\xi_{i+2}-\xi_{i}$, $i\in\M$,
yielding~(\ref{eq_at_most2}).

Next, iterating this argument, we obtain $|\xi_{i+2l}-\xi_i|\le
2l$. Because of the periodic boundary condition, in the case when
$M$ is odd, this results in $|\xi_i-\xi_j|\le 2M$ for all $i$ and
$j$, while in the case when $M$ is even this is true only whenever
$i-j$ is even. The observation that $\sum_j \xi_j(t)=t$ thus proves
(\ref{eq_xieo}) for odd $M$, since
$$
\left|t-M\xi_i(t)\right|=\left|\sum_{j=1}^{M}
\left[\xi_j(t)-\xi_i(t)\right]\right|\le (M-1)\times 2M.
$$
Now, when $M=2L$ is even, denote
 \bn
H(t)=\frac{\sum_{j=1}^L \xi_{2j}(t) - \sum_{j=1}^L
\xi_{2j-1}(t)}{M}.
 \en
Suppose $i\in\M$ is even. Then
 \bn
\left|t-M\xi_i(t)+MH(t)\right|&=&\left|2\sum_{j=1}^L
\xi_{2j}(t)-2L\xi_i(t) \right|\\
& =& \left|2\sum_{j=1}^L
\left[\xi_{2j}(t)-\xi_i(t)\right] \right| \le4M(L-1)<2M^2.
 \en
A similar argument holds for odd $i$. Hence we have
established~(\ref{eq_xieo}) for even $M$ as well as for odd $M$.

To finish the proof, denote by $\tau_m$, $m\ge 0$, the consecutive
renewal times of the Markov chain $v(t)$ after $t_{j_0}$, that is
 \bn
 \tau_0&=&\inf\{t\ge t_{j_0}:\ v_1(t)=v_2(t)=\dots=v_{M}(t)=0\},\\
 \tau_m&=&\inf\{t\ge \tau_{m-1}:\ v_1(t)=v_2(t)=\dots=v_{M}(t)=0\},\ m\ge 1.
 \en
By Proposition~\ref{prop_a_conv}, these stopping times are
well-defined; moreover, $\tau_{m+1}-\tau_m$ are i.i.d.\ and have
\index{tail behaviour}exponential tails. Let
$\zeta_{m+1}=H(\tau_{m+1})-H(\tau_{m})$.
Then the $\zeta_m$ are also i.i.d.; moreover, their distribution is
symmetric around $0$, and $|\zeta_{m+1}|\le \tau_{m+1}-\tau_m$
hence the $\zeta_m$ also have exponential tails. The rest follows
from the standard
\index{Donsker, M. D.!Donsker--Varadhan invariance principle}Donsker--Varadhan
invariance principle; see
\index{Durrett, R. T.}\index{Kesten, H.}\index{Limic, V.}e.g.~\cite{DKL},
pp.~590--592, for a proof in a very similar set-up.
\end{proof}

\section{Symmetric interaction}
In the symmetric case, the
\index{potential|)}potential of site $k$ at time $t$ is
 \bnn\label{eq:sympot}
  u_k(t)=\xi_{k-1}(t)+\xi_k(t)+\xi_{k+1}(t),\, k\in\M,
 \enn
and the transition probabilities of the Markov chain
$u(t)$ are now given by
 \begin{multline*}
 \P\left(u_{i}(t+1)=u_{i}(t)+1_{i\in\{k-1,k,k+1\}}, i=1,\ldots,M|u(t)\right) \\
=\left\{\begin{array}{rcl}
  0, &\mbox{if} & u_k(t)>m(t),\\ \\
  N^{-1}_{\min}(t), &\mbox{if} & u_k(t)=m(t),
 \end{array}\right.
 \end{multline*}
for $k\in\M$, where, as before,  $N_{\min}(t)\in\{1,2,\dots,M\}$ is the number
of $u_i(t)$ equal to $m(t)$.

\begin{proposition}\label{prop:by3}
If $(M\mbox{ mod }3)\ne 0$, then for any $u=(u_1,u_2,\dots,u_M)$
the system
 \bnn\label{eq:uxisym}
  u_1&=&\xi_M+\xi_{1}+\xi_{2}\nonumber\\
  u_2&=&\xi_1+\xi_{2}+\xi_{3}\nonumber\\
  &\vdots&\\
  u_M&=&\xi_{M-1}+\xi_{M}+\xi_{1}\nonumber
 \enn
has a unique solution. On the other hand, if $M$ is divisible by
$3$, system\/ \eqref{eq:uxisym} has a solution if and only if
 \begin{align}
 u_1+u_4+\dots+u_{M-2}&= u_2+u_5+\dots+u_{M-1}\nonumber\\
&= u_3+u_6+\dots+u_{M}.\label{eq:Mmod3}
 \end{align}
\end{proposition}

\begin{proof}If $M$ is not divisible by $3$, then the determinant
of the matrix
 \bn
 \left(
 \begin{array}{cccccccc}
 1 & 1 & 0 & 0 &\dots & 0 & 0 & 1 \\
 1 & 1 & 1 & 0 &\dots & 0 & 0 & 0 \\
 0 & 1 & 1 & 1 &\dots & 0 & 0 & 0 \\
 \vdots &\vdots &\vdots &\vdots &\ddots &\vdots &\vdots  &\vdots \\
 0 & 0 & 0 & 0 &\dots & 1 & 1 & 1 \\
 1 & 0 & 0 & 0 &\dots & 0 & 1 & 1 \\
 \end{array}
 \right)
 \en
corresponding to the equation~(\ref{eq:uxisym}) is $\pm 3$ (which
can be easily proved by induction). Hence the system must have a
unique solution.

When $M$ is divisible by $3$, by summing separately the
$1^\textup{st}$, $4^\textup{th}$, $5^\textup{th}$, \dots\ lines of
(\ref{eq:uxisym}), and then repeating this for the
$2\textup{nd}$, $5\textup{th}$, \dots\ or
$3\textup{th}$, $6\textup{th}$, \dots\ lines, we obtain
condition~(\ref{eq:Mmod3}).

Then it turns out that we can set both $\xi_1$ and $\xi_2$ to be
any real numbers, so $\xi_3=u_2-\xi_1-\xi_2$, and  $\xi_k$,
$k\ge 4$, are given:
 \bn
 \xi_{k+1}=[u_k-u_{k-1}]+[u_{k-3}-u_{k-4}]+\dots+\xi_{(k\mbox{ mod } 3)+1}.
 \en
\end{proof}

Similarly to the asymmetric case, consider the Markov chain $v(t)$  on
$\{0,1,2,\dots\}$ and  recall the definition of $t_j$
from~(\ref{eqv(t)}).
The following statement is straightforward.
\begin{proposition}\label{prop_sym1}
For any $k\in\M$
$$
v_{k}(t_{j+1})\le v_{k}(t_{j})
$$
unless both $v_{k-1}(t_{j})=0$ and $v_{k+1}(t_{j})=0$.
\end{proposition}

\begin{proposition}\label{prop_sym2}
For $j$ large enough, in the sequence of $v_k(t_j)$, $k\in\M$,
there are no more than two non-zero elements in a row, that is
$$
 \mbox{if }v_k(t_j)>0\mbox{ then either }v_{k-1}(t_j)=0\mbox{ or
 }v_{k+1}(t_j)=0, \mbox{ or both}.
$$
\end{proposition}

\begin{proof}Fix some $k\in\M$. Then  $v_k(t_j)$ is either $0$ or
positive. In the first case, unless both of the neighbours of
point $k$ are zeros at time $t_j$, by Proposition~\ref{prop_sym1}
we have $v_k(t_{j+1})=0$. On the other hand, if
$(v_{k-1}(t_j),v_k(t_j),v_{k+1}(t_j))=(0,0,0)$, then at time
$t_{j+1}$ either this triple becomes $(0,1,0)$ if both $k-1$ and
$k+1$ are chosen, or $v_k(t_{j+1})=0$.

Now suppose that $v_k(t_j)>0$. If both $v_{k-1}(t_{j})=0$ and
$v_{k+1}(t_{j})=0$, then from Proposition~\ref{prop_sym1} applied
to $k-1$ and $k+1$, we conclude
$v_{k-1}(t_{j+1})=v_{k+1}(t_{j+1})=0$, hence point $k$ remains
surrounded by $0$s. Similarly, if $v_k(t_j)>0$ and
$v_{k+1}(t_j)>0$ but $v_{k-1}(t_j)=v_{k+2}(t_j)=0$, then points
$\{k,k+1\}$ remain surrounded by $0$s at time $t_{j+1}$.

Finally, if point $k$ is surrounded by non-zeros on both sides,
that is $v_{k-1}(t_j)$, $v_{k}(t_j)$ and $v_{k+1}(t_j)$ are all
positive, we have $v_k(t_{j+1})=v_k(t_j)-1$.

Consequently, all sequences of non-zero elements of length $\ge 3$
are bound to disappear, and no such new sequence can arise as $j$
increases.
\end{proof}

\begin{proposition}\label{prop_sym3}
For any $k\in\M$, if for some $s$
$$
v_{k-1}(s)>0, \ v_{k}(s)=0, \ v_{k+1}(s)>0
$$
then for all $j$ such that $t_j\ge s$
$$
v_{k-1}(t_{j})>0, \ v_{k}(t_{j})=0, \ v_{k+1}(t_{j})>0.
$$
\end{proposition}

\begin{proof}This immediately follows from the fact that there
must be a particle
\index{adsorption|(}adsorbed at point $k$ during the time interval
$(t_{j_0},t_{j_0+1}]$ where $j_0=\max\{j:\ t_j\le s\}$, and that
would imply that $v_{k\pm 1}(t_{j_0+1})\ge v_{k\pm 1}(t_{j_0})$
while $v_{k}(t_{j_0+1})=0$. Now an induction on $j$ finishes the
proof.
\end{proof}

\begin{proposition}\label{prop_sym4}
For $j$ large enough, in the sequence of $v_k(t_j)$, $k\in\M$,
there are no more than two $0$s  in a row, that is
$$
 \mbox{if }v_k(t_j)=0\mbox{ then either }v_{k-1}(t_j)>0\mbox{ or
 }v_{k+1}(t_j)>0, \mbox{ or both}.
$$
\end{proposition}

\begin{proof}Suppose $j$ is so large that already there are no
consecutive subsequences of positive elements of length $\ge 2$ in
$(v_1,\dots,v_M)$ (see Proposition~\ref{prop_sym2}). Let
$$
Q(j)=\left|\left\{k:\ v_{k-1}(t_j)>0, v_{k}(t_j)=0,
v_{k-1}(t_j)>0, \right\}\right|.
$$
Proposition~\ref{prop_sym3} implies that $Q(j)$ is non-decreasing;
since $Q(j)<M$ it means that $Q(j)$ must converge to a finite
limit.

Let $A_j$ be the event that at time $t_j$ there are $3$ or more
zeroes in a row in $v(t_j)$. On $A_j$ there is a $k\in\M$ such
that $v_k(t_j)=v_{k+1}(t_j)=v_{k+2}(t_j)=0$ but $v_{k-1}(t_j)>0$,
(unless all $v_k=0$ but then the argument is similar).  Then, with
a probability exceeding $1/M$, at time $t_j+1$ new particle gets
adsorbed at $k+2$, yielding by Proposition~\ref{prop_sym3} that
for all $j'>j$ we have $v_{k-1}(t_{j'})>0$, $v_{k}(t_{j'})=0$,
$v_{k}(t_{j'})>0$, hence the event $B_j:=\{Q(j+1)\ge Q(j)+1\}$
occurs as well.  Therefore,
$$
\P(B_j\| \F_{t_j})\ge \frac 1M\, \P(A_j\| \F_{t_j}),
$$
where $\F_{t_j}$ denotes the sigma-algebra generated by $v(t)$ by
time $t_j$. Combining this with the second
\undex{Borel, F. E. J. E.!Borel--Cantelli lemma}Borel--Cantelli lemma,
we obtain
 \bn
 \{A_j \mbox{ i.o.}\}
  &=&\left\{\sum_j \P(A_j\| \F_{t_j})=\infty\right\}\subseteq
   \left\{\sum_j \P(B_j\| \F_{t_j})=\infty\right\}
  \\
  &=&\left\{B_j \mbox{ i.o.}\right\}
  =\left\{Q(j)\to\infty\right\}
 \en
leading to a contradiction.
\end{proof}

\begin{proposition}\label{prop_sym5}
Let
$$
W(j)=\left|\left\{k:\ v_{k-1}(t_j)=0, v_{k}(t_j)>0,
v_{k+1}(t_j)>0, v_{k-1}(t_j)=0, \right\}\right|
$$
be the number of `doubles'. Then $W(j)$ is non-increasing.
\end{proposition}

\begin{proof}
Let us investigate how we can obtain a subsequence
$(0,\ast,\ast,0)$ starting at position $k-1$ at time $t_{j+1}$,
where $\ast$ stands for a positive element. One possibility is
that at time $t_j$ we already have such a subsequence there; this
does not increase $W(j)$. The other possibilities at time $t_j$
are
 \bn
 (0,0,0,0),\  (0,0,0,\ast),\  (0,0,\ast,0 ),\ (0,0,\ast,\ast),\\
 (0,\ast,0,0),\  (0,\ast,0,\ast),\
 (\ast,\ast,0,0),\  (\ast,\ast,0,\ast).
 \en
By careful examination of all of the configurations above, we
conclude that the subsequence  $(0,\ast,\ast,0)$ cannot arise at
time $t_{j+1}$. Consequently, $W(j)$ cannot increase.
\end{proof}

\begin{proposition}\label{prop_sym6}
For $j$ large enough, in the sequence of $v_k(t_j)$, $k\in\M$,
there are no consecutive subsequences of the form
$(\ast,\ast,0,0)$ or $(0,0,\ast,\ast)$ where each $\ast$ stands for any
positive number; that is there is no $k$ such that
 \bn
v_{k}(t_j)=v_{k+1}(t_j)=0 \mbox{ and either }& & \\
&v_{k+2}(t_j)>0 \mbox{ and } v_{k+3}(t_j)>0& \\
\mbox{ or } &v_{k-1}(t_j)>0 \mbox{ and } v_{k-2}(t_j)>0.&
 \en
\end{proposition}

\begin{proof}Omitting $(t_j)$, without loss of generality suppose
$v_{k}>0$, $v_{k+1}>0$, $v_{k+2}=v_{k+3}=0$. Then either at some
time $j_1>j$ we will have  $v_{k+3}(t_{j_1})>0$ (hence the
configuration $(\ast,\ast,0,0)$ gets destroyed), or with
probability at least $\frac 13$ for each $j'\ge j$ we have
\index{adsorption|)}adsorption at position $k+3$ at some time during the time interval
$(t_{j'},t_{j'+1}]$. This would imply that
$v_{k+1}(t_{j'+1})=v_{k+1}(t_{j'})-1$. Hence, in a geometrically
distributed number of times, we obtain $0$ at position $k+1$, and
thus the configuration $(\ast,\ast,0,0)$ gets destroyed. On the
other hand, by Proposition~\ref{prop_sym5}, the number of doubles
is non-increasing, so no new configurations of this type can
arise. Consequently,  eventually all configurations
$(\ast,\ast,0,0)$ and $(0,0,\ast,\ast)$ will disappear.
\end{proof}

\begin{proposition}\label{prop_sym7}
For $j$ large enough, in the sequence of $v_k(t_j)$, $k\in\M$,
there are no consecutive subsequences of the form $(0,0,\ast,0,0)$
where $\ast$ stands for any positive number; that is there is no
$k$ such that
 \bn
v_{k-2}(t_j)=v_{k-1}(t_j)=0=v_{k+1}(t_j)=v_{k+2}(t_j) \mbox{ and }
v_{k}(t_j)>0.&
 \en
\end{proposition}

\begin{proof}Propositions~\ref{prop_sym2} and~\ref{prop_sym4}
imply that for some (random) $J$ large enough for all $j\ge J$
consecutive subsequences of zero (non-zero resp.) elements have
length $\le 2$, and Proposition~\ref{prop_sym6} says that two
consecutive $0$s must be followed (preceded resp.) by a single
non-zero element. Therefore, $(0,0,\ast,0,0)$ must be a part of a
longer subsequence of form $(0,\ast,0,0,\ast,0,0,\ast,0)$. This,
in turn, implies for the middle non-zero element located at $k$
that
 \bn
  v_k(t_{j+1})=\left\{\begin{array}{ll}
 v_k(t_j)+1, &\mbox{ with probability $1/4$},\\
 v_k(t_j), &\mbox{ with probability $1/2$},\\
 v_k(t_j)-1, &\mbox{ with probability $1/4$}.
  \end{array}\right.
 \en
Hence, by the properties of simple
\index{random walk (RW)}random walk, for some $j'>J$
we will have $v_k(t_{j'})=0$ (suppose that $j'$ is the first such
time). On the other hand, by Proposition~\ref{prop_sym1},
$$
v_{k-2}(t_{j'})=v_{k-1}(t_{j'})=v_{k+1}(t_{j'})=v_{k+2}(t_{j'})=0
$$
as well. This yields a contradiction with the choice of $J$ (see
Proposition~\ref{prop_sym4}).
\end{proof}

\begin{proof}[Proof of Theorem~\ref{thm:A3_b0}]Let a
\index{potential|(}\emph{configuration
of the potential} be a sequence $\bar v=(\bar v_1, \dots, \bar
v_M)$ where each $\bar v_i\in\{0,\ast\}$. Then we say that
$v=\break(v_1,v_2,\dots,v_M)$ with the following property has type $\bar v$:
 \bn
 v_i=0&\mbox{ if } &\bar v_i=0,\\
 v_i>0&\mbox{ if } &\bar v_i=\ast.
 \en

Propositions~\ref{prop_sym2}, \ref{prop_sym4}, \ref{prop_sym6},
and~\ref{prop_sym7} rule out various types of configurations for
all $j$ large enough. On the other hand, it is easy to check that
all remaining configurations for $v(t_j)$ are possible and stable,
that is, once you reach them, you stay in them forever.

Call a configuration $\bar v$ \emph{admissible}, if there is a
collection $\xi_1$, $\xi_2$, \dots, $\xi_M$ such that the
system~(\ref{eq:uxisym}) has a solution for some
$u=(u_1,\dots,u_M)$ having type $\bar v$. If $M$ is not divisible
by $3$, according to Proposition~\ref{prop_sym1} all
configurations $\bar v$ are admissible. On the other hand, it is
easy to see that if $M=3K$ then a necessary and sufficient
condition for a non-zero configuration $\bar v$ to be admissible
is
 \bn
 \bar v_i=\ast&\mbox{ for some $i$ such that }& i \mbox{ mod }3=0,\mbox{ and}\\
 \bar v_j=\ast&\mbox{ for some $j$ such that }& j \mbox{ mod }3=1,\mbox{ and}\\
 \bar v_k=\ast&\mbox{ for some $k$ such that }& k \mbox{ mod }3=2.
 \en
This establishes all possible stable configurations for $v$ and
hence the potential $u$,  thus determining the subset of $\M$
where points are adsorbed for sufficiently large times, namely,
$\xi_i(t)\to\infty$ if and only if $v_i(t_j)=0$ for all large $j$.

Moreover, whenever we see a subsequence of type
$(v_{k-1},v_k,v_{k+1})=(\ast,0,\ast)$, we have
$$
0\le \lim_{j\to\infty} \left[t_j-u_k(t_j)\right]<\infty,
$$
and for a subsequence of type
$(v_{k-1},v_k,v_{k+1},v_{k+2})=(\ast,0,0,\ast)$ we have
$$
\lim_{j\to\infty} \frac{u_k(t_j)}{t_j}
 =\lim_{j\to\infty} \frac{u_{k+1}(t_j)}{t_j}=\frac 12
$$
by the strong law. Setting
$$
\a=\frac 1{\lim_{j\to\infty} |\{i\in\M:\ v_i(t_j)>0,
v_{i+1}(t_j)=0\}|}
$$
finishes the proof of the first part of the Theorem.

Finally, note that if the initial configuration is empty, the
conditions of Proposition~\ref{prop_sym5} are fulfilled with no
`doubles' at all, i.e.\ $W(0)=0$. Consequently, for all $j\ge 0$
we have that there are no consecutive non-zero elements in
$v_k(t_j)$, yielding the final statement of the Theorem.
\end{proof}

\section{Appendix}\label{binf}\index{allocation process|(}

 In this section we briefly describe the long-time  behaviour
of the growth process  generated by the dynamics, where a particle is
allocated at random to a site with maximum potential. The
process is trivial in both the symmetric and asymmetric cases. Consider
the symmetric case, i.e.  $U_i=\{i-1,i,i+1\}$, $i\in\M$. It is
easy  to see that with probability $1$, there exists $k$ such that
either
 \begin{equation}
 \label{sym_inf}
\lim_{t\to\infty} \frac{\xi_k(t)}t=1 \ \mbox{ and } \ \sup_{i\ne k}
 \xi_i(t)<\infty
 \end{equation}
or
 \begin{equation}
\label{sym_inf1}
 \lim_{t\to\infty} \frac{\xi_k(t)}t=\lim_{t\to\infty} \frac{\xi_{k+1}(t)}t=\frac 12
  \ \mbox{ and } \ \sup_{i\notin \{k,k+1\}}  \xi_i(t)<\infty.
 \end{equation}
Indeed, recall the formula for the
\index{potential|)}potential given by~(\ref{eq:sympot}).
Then $u(t)$ is a
\index{Markov, A. A.!Markov chain|)}Markov chain with transition
probabilities given by
 \bn
 \P\left(u_{i}(t+1)=u_{i}(t)+1_{i\in\{k-1,k,k+1\}}\right)=\frac{1_{\{k\in S_{\max}(t)\}}} {|S_{\max}(t)|}
 \en
for $k\in\M$, where
$$
S_{\max}=\left\{i:\ u_i(t)=\max_{i\in\M} u_i(t)\right\}\subseteq
\M
$$
is the set of those $i$ for which $u_i(t)$ equals the maximum value.

Observe that if at time $s$ the
\index{adsorption}adsorption/allocation  occurs at point $i$,
then $S_{\max}(s+1)\subseteq \{i-1,i,i+1\}$. In particular, if the
maximum is unique, that is, $S_{\max}(s+1)=\{i\}$, then for all
times $t\ge s$ this property will hold, and hence all the
particles from now on will be adsorbed at $i$ only.

If, on the other hand, $|S_{\max}(s+1)|=2$, without loss of
generality say $S_{\max}(s+1)=\{i,i+1\}$, then  this property
will be also  preserved for all $t>s$ and each new particle will be
adsorbed with probability $\frac 12$ at either $i$ or $i+1$.

Finally, if $|S_{\max}(s+1)|=3$, say
$S_{\max}(s+1)=\{i,i+1,i+2\}$, then at time $s+2$ either
$S_{\max}(s+2)=\{i,i+1,i+2\}$ if the adsorption occurred at $i+1$,
or $S_{\max}(s+2)=\{i+1,i+2\}$ or $\{i,i+1\}$ otherwise. By
iterating this argument we obtain that after a geometric number of
times we will arrive at the situation where $|S_{\max}(t)|=2$, and
then the process will follow the pattern described in the previous
paragraph.

A similar simple argument shows that in the case of the asymmetric
interaction only the outcome (\ref{sym_inf}) is
possible.\index{allocation process|)}

\end{document}